\documentclass[11pt]{amsart}
\input{diagrams}

\newtheorem{proposition}{Proposition}[section]
\newtheorem{lemma}[proposition]{Lemma}
\newtheorem{corollary}[proposition]{Corollary}
\newtheorem{theorem}[proposition]{Theorem}

\theoremstyle{definition}

\theoremstyle{remark}
\newtheorem{remark}[proposition]{Remark}

\newcommand{\thlabel}[1]{\label{th:#1}}
\newcommand{\thref}[1]{Theorem~\ref{th:#1}}
\newcommand{\selabel}[1]{\label{se:#1}}
\newcommand{\seref}[1]{Section~\ref{se:#1}}
\newcommand{\lelabel}[1]{\label{le:#1}}
\newcommand{\leref}[1]{Lemma~\ref{le:#1}}
\newcommand{\prlabel}[1]{\label{pr:#1}}
\newcommand{\prref}[1]{Proposition~\ref{pr:#1}}
\newcommand{\colabel}[1]{\label{co:#1}}
\newcommand{\coref}[1]{Corollary~\ref{co:#1}}
\newcommand{\relabel}[1]{\label{re:#1}}
\newcommand{\reref}[1]{Remark~\ref{re:#1}}

\newcommand{\eqlabel}[1]{\label{eq:#1}}
\newcommand{\equref}[1]{(\ref{eq:#1})}

\newcommand{\Hom}{{\rm Hom}}
\newcommand{\HOM}{{\rm HOM}}
\newcommand{\End}{{\rm End}}
\newcommand{\Ext}{{\rm Ext}}
\newcommand{\Tor}{{\rm Tor}}
\newcommand{\EXT}{{\rm EXT}}

\newcommand{\Ker}{{\rm Ker}\,}

\newcommand{\im}{{\rm Im}\,}

\newcommand{\Spec}{{\rm Spec}\,}

\def\ot{\otimes}

\newcommand{\Mm}{\mathcal{M}}

\def\text#1{{\rm {\rm #1}}}

\def\ol{\overline}

\begin{document}
\title[On the cohomology of $H$-comodules]{On the cohomology of relative Hopf
modules}
\author{S. Caenepeel}
\address{Faculty of Applied Sciences,
Vrije Universiteit Brussel, VUB, B-1050 Brussels, Belgium}
\email{scaenepe@vub.ac.be}
\urladdr{http://homepages.vub.ac.be/\~{}scaenepe/}
\author{T. Gu\'ed\'enon}
\address{Faculty of Applied Sciences,
Vrije Universiteit Brussel, VUB, B-1050 Brussels, Belgium}
\email{tguedeno@vub.ac.be}
\thanks{Research supported by the project G.0278.01 ``Construction
and applications of non-commutative geometry: from algebra to physics"
from FWO Vlaanderen}
\begin{abstract}
Let $H$ be a Hopf algebra over a field $k$, and $A$ an $H$-comodule algebra.
The categories of comodules and relative Hopf modules are then Grothendieck categories
with enough injectives.
We study the derived functors of the associated Hom functors, and of the
coinvariants functor, and discuss spectral sequences that connect them. We
also discuss when the coinvariants functor preserves injectives.
\end{abstract}

\maketitle

\section*{Introduction}
Let $k$ be a field, and $H$ a Hopf algebra with bijective antipode,
and $A$ an $H$-module algebra. We can then consider the smash product
$A\# H$ and the subring of invariants $A^H$. A  left $H$-module $M$
is called locally finite if $\dim_k(Hm)$ is finite, for every $m\in M$.
In \cite{18}, the second author studied homological algebra for $H$-locally
finite $A\#H$-modules, with emphasis to injective modules, minimal injective
resolutions and cohomology. He also
calculated the Picard group of $A^H$ in terms of the Picard group of $A$
and various subgroups of the group
$Z(H , A)$ consisting of linear maps from $H\to A$ satisfying the cocycle
condition. In the particular situation where $H$ is the enveloping algebra
of a finite dimensional Lie algebra, we refer to \cite{16,17}. The methods
in \cite{16,17,18} are based on Magid's papers \cite{22,23} on rational algebraic
group actions.

The aim of this paper is to discuss the homological algebra for relative Hopf modules.
If $H$ is a Hopf algebra, and $A$ is an $H$-comodule algebra, then a relative
Hopf module is a vector space with an $A$-action and an $H$-coaction with a certain
compatibility relation. In the case where $H$ is finite dimensional, the category of
relative Hopf modules is isomorphic to the category of modules over the smash
product $A\# H^*$, providing the connection to the theory developed in \cite{18}.
However, the situation is more interesting in the case where $H$ is infinite dimensional.
Given two $H$-comodules $M$ and $N$, we can consider the space $\Hom^H(M,N)$
of $H$-colinear morphisms between $M$ and $N$, and also the $H$-comodule
$\HOM(M,N)$, consisting of rational $k$-linear maps $M\to N$. We can consider
the right derived functors of these two Hom functors, given rise to two different versions
of the Ext functors. $\Hom^H$ can be viewed as the composition of $\HOM$ and the
coinvariants functor, and the results in a spectral sequence connecting the two versions
of Ext. This is discussed in \seref{1}. In \seref{2}, we look at relative Hopf modules. Again,
we have two versions of the Hom functor, and, with some addtional conditions, the corresponding right derived functors
are connected by a spectral sequence, see Propositions \ref{pr:2.11} and \ref{pr:2.15}.
More specific results can be obtained in the case where $H$ is cosemisimple, this is
discussed in \seref{3}.

\section{The right derived functors of the coinvariant functor and the HOM 
functor}\selabel{1}
Throughout this paper, $k$ is a field, and $H$ is a Hopf algebra with bijective
antipode. We recall that $\Mm^H$, the category of $H$-comodules and $H$-colinear
maps, is a Grothendieck category with enough injectives (see for example \cite{11}).
We say that $H$ has the symmetry property if $M\ot N$ and $N\ot M$ are isomorphic
as comodules, for any $M,N\in \Mm^H$. If $H$ is an almost commutative Hopf algebra,
then its antipode is bijective and it has the symmetry property (see
\cite[10.2.11,~10.2.12]{24}).\\
Let $A$ be an $H$-comodule algebra. A relative left-right $(A,H)$-Hopf module
is a vector space with a left $A$-action and a right $H$-coaction $\rho$ such that
$\rho(am)=\rho(a)\rho(m)$, for all $a\in A$ and $m\in M$. The category of relative
$(A,H)$-Hopf modules ${}_A\Mm^H$ has direct sums, and is a Grothendieck category with
enough injective objects. If $A$ is noetherian, then direct sums of injectives
are injective, see \cite[3.1, 3.2]{29}.\\
We will use the Sweedler-Heyneman notation for comultiplications and coactions;
if $\Delta$ is the comultiplication on $H$, then we write
$$\Delta(h)=h_1\ot h_2,$$
where the summation is implicitly understood. In a similar way, if $M$ is
a right $H$-comodule, with right $H$-coaction $\rho$, then we write, for all
$m\in M$:
$$\rho(m)=m_0\ot m_1.$$
$M^{{\rm co}H}=\{m\in M~|~\rho(m)=m\ot 1\}$ is called the $k$-submodule of
coinvariants of $M$. $\Mm^H$ is a monoidal category: if $M,N\in \Mm^H$,
then $M\ot N\in \Mm^H$, with $H$-coaction
$$\rho(m\ot n)=m_0\ot n_0\ot m_1n_1.$$
The unit object is $k$, with coaction $\rho(x)=x\ot 1_H$.\\
Take $f\in \Hom(M,N)$, and consider $\rho(f)\in \Hom(M,N\ot H)$ given by
$$\rho(f)(m)=f(m_0)_0\ot S^{-1}(m_1)f(m_0)_1.$$
As $k$ is a field, $\Hom(M,N)\ot H\subset \Hom(M,N\ot H)$, and we introduce
$$\HOM(M,N)=\{f\in \Hom(M,N)~|~\rho(f)\in \Hom(M,N)\ot H\}.$$
A morphism $f\in \HOM(M,N)$ is called a rational morphism. 
If $H$ is finite dimensional, then all morphisms are rational. It is well-known
(see for example \cite{25,28}) that $\HOM(M,N)$ is an $H$-comodule, and that
it is the largest $H$-comodule contained in $\Hom(M,N)$. Also recall that
$\rho(f)=f_0\ot f_1$ if and only if
\begin{equation}\eqlabel{1.1.1}
f_0(m)\ot f_1= f(m_0)_0\ot S^{-1}(m_1)f(m_0)_1.
\end{equation}

\begin{lemma}\lelabel{1.1}
For any $M,N\in \Mm^H$, we have that $\HOM(M,N)^{{\rm co}H}=\Hom^H(M,N)$.
\end{lemma}

\begin{proof}
If $\rho(f)=f\ot 1$, then it follows from \equref{1.1.1} that
$$f(m_0)\ot m_1=f(m_0)_0\ot m_2S^{-1}(m_1)f(m_0)_1=f(m)_0\ot f(m)_1$$
and $f$ is $H$-colinear. Conversely, if $f$ is $H$-colinear, then
$$f(m_0)_0\ot S^{-1}(m_1)f(m_0)_1=f(m_0)\ot S^{-1}(m_2)m_1=f(m)\ot 1,$$
and it follows from \equref{1.1.1} that $\rho(f)=f\ot 1$.
\end{proof}

\begin{proposition}\prlabel{1.2}
Let $M,N,P\in\Mm^H$, and consider the
natural isomorphism of vector spaces
$$\phi:\ \Hom(N\otimes M, P)\to Hom(M, Hom(N, P)),~~
\phi(f)(m)(n)=f(n\otimes m).$$
\begin{enumerate}
\item If $f\in \Hom(N\otimes M, P)$ is $H$-colinear, then 
$\phi(f)(m)\in \HOM(N, P)$, for every $m \in M$; furthermore $\phi(f)$ is $H$-colinear.
\item $\phi$ induces a $k$-isomorphism
$$ \phi:\ \Hom^H(N\otimes M, P)\to \Hom^H(M, \HOM(N,  P)).$$
\item If $H$ has the symmetry property, then $\phi$ induces a $k$-isomorphism
$$ \psi:\ \Hom^H(M\otimes N, P)\to \Hom^H(M, \HOM(N,  P)).$$
\end{enumerate}
\end{proposition}

\begin{proof}
(1) Let $f$ be $H$-colinear. We claim that
$$\rho(\phi(f)(m))=f(-\ot m_0)\ot m_1.$$
Indeed, we show easily that \equref{1.1.1} is satisfied:
\begin{eqnarray*}
&&\hspace*{-2cm}
(\phi(f)(m)(n_0))_0
\otimes S^{-1}(n_1)(\phi(f)(m)(n_0))_1\\
&=&f(n_0\otimes m)_0 \otimes S^{-1}(n_1)f(n_0 \otimes m)_1\\
&=&f(n_{00} \otimes m_0) \otimes S^{-1}(n_1)n_{01}m_1\\
&=&f(n_0\otimes m_0)
\otimes S^{-1}(n_2)n_1m_1=f(n\otimes m_0) \otimes m_1,
\end{eqnarray*}
as needed. $\phi(f)$ is $H$-colinear
$$\phi(f)(m_0)\ot m_1=\rho(\phi(f)(m))=f(-\ot m_0)\ot m_1,$$
for all $m\in M$. This is equivalent to
$$\phi(f)(m_0)(n)\ot m_1=f(n\ot m_0)\ot m_1,$$
for all $m\in M$ and $n\in N$. This is obvious.\\

(2) Take $f:N\ot M\to P$, and assume that
$\phi(f)\in \Hom^H(M, \HOM(N,  P))$. Then we compute that
\begin{eqnarray*}
&&\hspace*{-2cm}
\rho(f(n\ot m))=\rho((\phi(f)(m))(n)=\rho\Bigl((\phi(f)(m))(n_0)\Bigr)\varepsilon(n_1)\\
&=&\Bigl((\phi(f)(m))(n_0)\Bigr)_0\ot n_2S^{-1}(n_1)   \Bigl((\phi(f)(m))(n_0)\Bigr)_1\\
&=& f(n_0\ot m_0)\ot n_1m_1,
\end{eqnarray*}
and it follows that $f$ is right $H$-colinear.\\

By the symmetry property, there is an $H$-colinear isomorphism
$\tau:\ N\ot M\to M\ot N$. The map
$$\Hom^H(\tau,P):\ \Hom^H(M\ot N,P)\to \Hom^H(N\ot M,P)$$
is an isomorphism of vector spaces, and $\psi= \phi\circ \Hom^H(\tau,P)$
is the required isomorphism.
\end{proof}

\begin{corollary}\colabel{1.3}
Let $M,V\in \Mm^H$, with $V$ finite dimensional, and $M$ projective in $\Mm^H$.
Then $V\ot M\in \Mm^H$ is also projective.
\end{corollary}

\begin{proof}
As $V$ is finite dimensional, the $H$-comodules $\HOM(V,P)=\Hom(V,P)$
and $V^*\ot P$ are isomorphic, for all $P\in \Mm^H$. Therefore
$\Hom(V,-):\ \Mm^H\to \Mm^H$ is exact. Also $\Hom^H(M,-):\ \Mm^H\to \Mm$
is exact, since $M\in \Mm^H$ is projective. It then follows from \prref{1.2} (2)
that $\Hom^H(V\ot M,-):\ \Mm^H\to \Mm^H$ is exact, and $V\ot M$ is a projective
object in $\Mm^H$.
\end{proof}

Recall that $M\in \Mm^H$ is called simple if it has no proper subobjects.
$M$ is semisimple or completely reducible if it is isomorphic to the direct sum
of simple objects. $\Mm^H$ is called semisimple or completely reducible if
every object is semisimple. It is well-known that
$\Mm^H$ is semisimple if and only if $H$ is
cosemisimple, see for example \cite[Lemma 2.4.3]{24}. We present another criterion
in the \leref{1.4}.

\begin{lemma}\lelabel{1.4}
$\Mm^H$ is semisimple if and only if $k$ is a projective object in
$\Mm^{{\rm fd}H}$, the category
of finite dimensional $H$-comodules.
\end{lemma}

\begin{proof}
Take an exact sequence
$$0 \rightarrow W_1 \rightarrow W_2 \rightarrow W_3 \rightarrow 0$$
in $\Mm^{{\rm fd}H}$, and let $V$ be a finite dimensional right $H$-comodule.
Then we have the following exact sequence in $\Mm^{{\rm fd}H}$:
$$0\rightarrow \Hom(V , W_1)\rightarrow \Hom(V , W_2)\rightarrow \Hom(V , W_3)
\rightarrow 0.$$
If $k$ is a projective object in $\Mm^{{\rm fd}H}$, then we have the following
exact sequence of vector spaces
$$\begin{matrix}
0&\rightarrow &\Hom^H(k , \Hom(V , W_1))&\rightarrow& \Hom^H(k , \Hom(V , W_2))\\
&\rightarrow& \Hom^H(k , \Hom(V ,W_3))&\rightarrow& 0
\end{matrix}$$
It follows from (2) in \prref{1.2} that the sequence
$$0\rightarrow \Hom^H(V , W_1)\rightarrow \Hom^H(V , W_2)\rightarrow  \Hom^H(V
, W_3))\rightarrow 0$$
is exact, so $V$ is a projective object in $\Mm^{{\rm fd}H}$, and
therefore any
subcomodule of $V$ is a direct summand of $V$ in $\Mm^H$. It follows
that $V$ is semisimple in $\Mm^H$.
Let $M$ be in $\Mm^H$. By the Fundamental Theorem of comodules 
\cite[Theorem 2.1.7]{11}, each element $m\in
M$ is contained in a finite-dimensional subcomodule $V_m$ of $M$. In
particular, every $m\in M$ is
contained in a sum of simple subcomodules of $M$, this implies that $M$ is the sum of
a family of simple subojects.
Using Zorn's Lemma we can show that this sum is direct.
\end{proof}

Using \prref{1.2}, we now give necessary and sufficient conditions for
the rationality of $f\in \Hom(M,N)$.

\begin{proposition}\prlabel{1.5}
Take two $H$-comodules $M$ and $N$.
For $f\in \Hom(M,N)$, the following assertions are equivalent.
\begin{enumerate}
\item $f\in \HOM(M,N)$;
\item there exists an $H$-comodule $V$, an element $v$ in
$V$ and an $H$-colinear map $F:\ M\otimes V\to N $ such that
$F(m\otimes v)=f(m)$ for all $m$ in $M$;
\end{enumerate}
If $H$ has the symmetry property, then (1) and (2) are equivalent to
\begin{enumerate}
\item[(3)] there exists an $H$-comodule $V$, an element $v$ in
$V$ and an $H$-colinear map $F':\ M\to \Hom(V, N)$ such that
$F'(m)(v)=f(m)$ for all $m$ in $M$.
\end{enumerate}
In (2) and (3), we can choose $V$ to be finite dimensional.
\end{proposition}

\begin{proof}
(2)$\Longleftrightarrow$(3) follows from \prref{1.2}.\\

(2)$\Longrightarrow$(1). 
We claim that $\rho(f)=F(-\ot v_0)\ot v_1$. Using the $H$-colinearity of $F$,
we obtain
\begin{eqnarray*}
&&\hspace*{-2cm}
f(m_0)_0\ot S^{-1}(m_1)f(m_0)_1=F(m_0\ot v)_0\ot S^{-1}(m_1)F(m_0\ot v)_1\\
&=& F(m_0\ot v_0)\ot S^{-1}(m_2)m_1v_1=F(m\ot v_0)\ot v_1
\end{eqnarray*}
and \equref{1.1.1} holds, as needed.\\

(1)$\Longrightarrow$(2). Take a finite dimensional $H$-subcomodule of
$\HOM(M,N)$ containing $f$. Such a $V$ exists by 
the Fundamental Theorem \cite[Theorem 2.1.7]{11}.
Then define $F:\ M\ot V\to N$ by
$$F(m\ot v)=v(m)$$
Clearly $F(m\ot f)=f(m)$, so we are done if we can show that $F$ is $H$-colinear.
Using the fact that $v\in V$ is rational, we find
\begin{eqnarray*}
&&\hspace*{-2cm} F(m_0\ot v_0)\ot m_1v_1=v_0(m_0)\ot m_1v_1\\
&=& v(m_0)_0\ot m_2S^{-1}(m_1)v(m_0)_1= v(m)_0\ot v(m)_1\\
&=& F(m\ot v)_0 \ot  F(m\ot v)_1.
\end{eqnarray*}
\end{proof}

\begin{corollary}\colabel{1.6}
Take $M, N,P\in \Mm^H$ be $H$-comodules. If $g\in \HOM(M, N)$ and
$f\in \HOM(N, P)$, then $f\circ
g\in \HOM(M, P)$.
\end{corollary}

\begin{proof} 
By \prref{1.5}, there exist finite dimensional $H$-comodules $V$ and
$W$, $v\in V$, $w \in
W$ and $H$-colinear maps $G:\ M \otimes V \rightarrow N$, $F:\ N
\otimes W \rightarrow
P$ such that $G(m \otimes v)=g(m)$,  $F(n \otimes w)=f(n)$ for all $m \in
M$ and $n \in N$. The map
$$K:\ M \otimes V \otimes W \rightarrow P,~~K(m
\otimes s \otimes t)= F(G(m
\otimes s) \otimes t)$$
is $H$-colinear, and $K(m\ot(v\ot w))=(f\circ g)(m)$.
\end{proof}

\begin{corollary}\colabel{1.7}
For any $T\in \Mm^H$,
$\HOM(-,T)$ and $\HOM(T,-)$ are left exact endofunctors of $\Mm^H$.
\end{corollary}

\begin{proof} 
Let 
$0 \rightarrow M \rTo^{i} N \rTo^{\pi} P \rightarrow
0$ be an exact sequence
in $\Mm^H$. Then
 $$0 \rightarrow \Hom(P, T) \rightarrow \Hom(N, T)
\rightarrow \Hom(M, T)\rightarrow
0$$ 
is an exact sequence of vectorspaces. $\pi$ is $H$-colinear, hence $\pi\in \HOM(N, P)$,
by \leref{1.1}. Consequently, $f\circ\pi\in \HOM(N, T)$, for all $f\in \HOM(P,T)$.
In a similar way, $f \circ i \in \HOM(M, T)$, for all  $f
\in \HOM(N, T)$, and it follows that $\HOM(-, T)$ is left exact.
\end{proof}

\begin{proposition}\prlabel{1.8}
Let $I$ be an injective object of $\Mm^H$. Then
\begin{enumerate}
\item $\HOM(N, I)$ is an injective object of $\Mm^H$, for any $N\in \Mm^H$;
\item $\HOM(-, I)$ is an exact endofunctor of $\Mm^H$.
\end{enumerate}
\end{proposition}

\begin{proof}
(1) follows from \prref{1.2} and the fact that $N\ot -:\ \Mm^H\to \Mm^H$ is
exact.\\
(2) Let 
$0 \rightarrow M \rTo^{i} N \rTo^{\pi} P \rightarrow
0$ be an exact sequence in $\Mm^H$. We know from \coref{1.7} that
$$0 \rightarrow \HOM(P, I) \rightarrow \HOM(N, I)
\rightarrow \HOM(M, I)$$ 
is exact in $\Mm^H$. Take $f\in \HOM(M, I)$, and let $V$ be a finite dimensional
$H$-subcomodule of $\HOM(M, I)$ containing $f$. Clearly
$i\ot V:\ M\ot V\to N\ot V$ is an $H$-colinear monomorphism. As in the proof
of $(1)\Rightarrow (2)$ in \prref{1.5}, we can show that
$$F:\ M\ot V\to I,~~F(m\ot v)=v(m)$$
is rational. Since $I\in \Mm^H$ is injective, there exists an $H$-colinear
map $G:\ N\ot V\to I$ such that $G\circ (i\ot V)=F$. It follows from \prref{1.5}
that
$$g:\ N\to I,~~g(n)=G(n\ot f)$$
is rational. On the other hand
$$f(m)=F(m\ot f)=G((i\ot V)(m\ot f))=G(i(m)\ot f)=(g\circ i)(m),$$
and it follows that $\HOM(N, I)
\rightarrow \HOM(M, I)$ is surjective.
\end{proof}

We will use the following notation.
\begin{itemize}
\item $R^pa({\rm co}H,-)$ are the right derived functors of the covariant left
exact functor $(-)^{{\rm co}H}:\ \Mm^H\to \Mm$;
\item $\EXT^p(M,-)$ are the right derived functors of $\HOM(M,-):\ \Mm^H\to \Mm^H$;
\item $\Ext^{H^p}(-,-)$ are the right derived functors of
$\Hom^H(-,-):\ \Mm^H\times \Mm^H\to \Mm$.
\end{itemize}
In particular, if $M$ and $N$ are $H$-comodules, then $\EXT^p(M,N)$ is also an
$H$-comodule. If $V\in \Mm^H$ is finite dimensional, then
$\Hom(V,-)\cong V^*\ot M$, hence $\HOM(V,-)$ is exact, and
$\EXT^q(V,M)=0$ for all $q\geq 1$.

\begin{proposition}\prlabel{1.9}
Let $M,N\in \Mm^H$.
\begin{enumerate}
\item We have a spectral sequence
$$
R^pa({\rm co}H, \EXT^q(M,N))~~\Rightarrow~~\Ext^{H^{p+q}}(M,N)
$$
with $p,q\geq 0$;
\item if $M$ is finite dimensional, then
$$
R^pa({\rm co}H, M^*\ot N)=\Ext^{H^p}(M,N),
$$
for all $p\geq 0$.
\end{enumerate}
\end{proposition}

\begin{proof}
By \leref{1.1}, $\HOM(M,N)^{{\rm co}H}=\Hom^H(M,N)$, and the result follows from
\prref{1.8} (1) and Grothendieck's spectral sequence for composite functors.
\end{proof}

In order to be able to compute right derived functors, we describe injective
resolutions of $M\in \Mm^H$. 

Let $V$ be a vector space. Then $V\ot H$ is a right $H$-comodule, with coaction
induced by the comultiplication, and we call $V\ot H$ a free $H$-comodule.
Recall \cite[Prop. 2.4.7]{11} that a right $H$-comodule $M$ is an injective object
in $\Mm^C$ if and only if it is a direct summand in a free $H$-comodule. In particular
$H$ is injective. \leref{1.10} is the analog of
\cite[Prop. 3.10 (c)]{19} for the category of $H$-comodules.

\begin{lemma}\lelabel{1.10}
Take $M,N\in \Mm^H$.
If $N\in \Mm^H$ is injective, then $M\ot N$ is also injective.
In particular, $M\ot H$ is an injective object of $\Mm^H$.
\end{lemma}

\begin{proof}
As we have seen above, $N$ is a direct summand of $V\ot H$, with $V$ a vector
space. Then $M\ot N$ is a direct summand of $M\ot V\ot H$.
Let $M_{\rm tr}$ be the vector space $M$ with trivial $H$-coaction.
We have an isomorphism of $H$-comodules
$$f:\ M\ot V\ot H\to M_{\rm tr}\ot V\ot H,~~f(m\ot v\ot h)=m_{0}\ot v\ot m_{1}h,$$
with inverse given by $f^{-1}(m\ot v\ot h)=m_{0}\ot v\ot S(m_{1})h$. So
$M\ot N$ is a direct summand of the free comodule $M_{\rm tr}\ot V\ot H$,
and is an injective object of $\Mm^H$.
\end{proof}

For $M\in \Mm^H$, we define $C^q(M)$ and $\varphi_q:\ C^q(M)\to C^{q+1}(M)$ recursively by
$$C^{-1}(M)=M~~{\rm and}~~C^{q+1}(M)=C^q(M)\ot H;$$
$$\varphi_{-1}:\ M\to M\ot H,~~\varphi_{-1}(m)=m\ot 1;$$
$$\varphi_{q+1}(u\ot h)=u\ot h\ot 1-\varphi_q(u)\ot h.$$
It is clear that $\varphi_q$ is $H$-colinear. Using induction on $q$, we easily
show that $\varphi_{q+1}\circ \varphi_q=0$, hence $\{C^q(M)\}_{q\geq 0}$ is
a complex in $\Mm^H$. Now consider
$$\psi_q:\ C^q(M)\to C^{q-1}(M),~~\psi_q(u\ot h)=\varepsilon(h)u.$$
Then a straightforward computation shows that
$$\varphi_{q-1}\circ\psi_q+\psi_{q+1}\circ\varphi_q=C^q(M),$$
the identity map on $C^q(M)$, for all $q\geq 0$. Hence
$\im (\varphi_q)\supset \Ker(\varphi_{q+1})$, and $C^*(M)$ is an acyclic complex.
It follows from \leref{1.10} that $C^q(M)$ is an injective object in $\Mm^H$,
for all $q\geq 0$, hence $C^*(M)$ is an injective resolution of $M\in \Mm^H$.
It follows that $R^pa({\rm co}H, M)$ is the cohomology group of the complex
${C^{*}(M)}^{{\rm co}H}$, and $\EXT^{p}(M, N)$ is
the cohomology group of the complex
$\HOM(M, C^*(N))$.

\section{The right derived functors of ${}_A\HOM(-, -)$ and ${}_A\Hom^H(-, -)$}\selabel{2}
Let $A$ be a right $H$-comodule algebra. Recall that this is an algebra
with a right $H$-coaction $\rho_A$ such that the unit and the multiplication are
right $H$-colinear, that is,
$${\rho}_A (ab) = a_0b_0 \otimes a_1b_1 \quad \hbox {and} \quad {\rho}_A
(1_A) = 1_A \otimes 1_H.$$
A vector space $M$ with a left $A$-action and a right $H$-coaction $\rho_M$ is called
a relative $(A,H)$-Hopf module if
$${\rho}_M(am)=a_0m_0 \otimes a_1m_1,$$
for all $a\in A$ and $m\in M$. ${}_A\Mm^H$ is the category of relative Hopf module
and $A$-linear $H$-colinear maps. For two relative Hopf modules $M$ and $N$,
we let ${}_A\Hom^H(M,N)$ be the space of $A$-linear $H$-colinear maps, and
$${}_A\HOM(M, N)={}_A\Hom(M, N)\cap \HOM(M, N).$$
The aim of this Section is to relate the right derived functors of
${}_A\HOM(-, -)$ and ${}_A\Hom^H(-, -)$ by a spectral sequence.
The sequence collapses if $H$ is cosemisimple. We can improve the results if
$A$ is left noetherian.

\begin{lemma}\lelabel{2.1}
Let $M$ and $N$ be relative $(A, H)$-Hopf modules, and take $f\in {}_A\Hom(M, N)$.
\begin{enumerate}
\item The $k$-linear map $\rho(f) :\ M\rightarrow
N\otimes H$ defined by
$$\rho(f)(m)=f(m_0)_0 \otimes S^{-1}(m_1)f(m_0)_1$$
is $A$-linear, hence $\rho(f)\in {}_A\Hom(M, N\otimes H)$;
\item $f\in {}_A\HOM(M , N)$ if and only if there exists  $f_0 \otimes
f_1\in {}_A\Hom(M , N)\otimes H$
such that $$f_0(m) \otimes f_1=f(m_0)_0 \otimes S^{-1}(m_1)f(m_0)_1,$$
for all $m\in M$.
\end{enumerate}
\end{lemma}

\begin{proof}
For all $a\in A$ and $m\in M$, we have
\begin{eqnarray*}
\rho(f)(am)&=& f(a_0m_0)_0 \otimes S^{-1}(a_1m_1)f(a_0m_0)_1\\
&=&a_0f(m_0)_0\ot S^{-1}(m_1)S^{-1}(a_2)a_1f(m_0)_1\\
&=&a(\rho(f)(m)).
\end{eqnarray*}
This shows that $\rho(f)\in {}_A\Hom(M, N\otimes H)$. The second part then
follows easily.
\end{proof}

\begin{lemma}\lelabel{2.2}
Let $M$ and $N$ be relative $(A , H)$-Hopf modules.
\begin{enumerate}
\item ${}_A\HOM(M, N)$ is an $H$-subcomodule of $\HOM(M, N)$; it is the
largest $H$-comodule contained in
${}_A\Hom(M,  N)$;
\item ${}_A\HOM(M,  N)^{{\rm co}H}={}_A\Hom^H(M, N)$;
\item if $M\in {}_A\Mm$ is finitely generated, then ${}_A\HOM(M,  N)={}_A\Hom(M,
N)$.
\end{enumerate}
\end{lemma}

\begin{proof} 
(1) follows from \leref{2.1} and the fact that $\HOM(M,N)$ is the largest comodule
contained in $\Hom(M,N)$ (see \seref{1}).\\
(2) similar to the proof of \leref{1.1}.\\ 
(3) follows from \cite[Prop. 4.2]{9}.
\end{proof}

Now let $M\in \Mm^H$, and $N\in {}_A\Mm^H$. We have seen in \cite[Lemma 1.1]{9}
that $N\ot M\in {}_A\Mm^H$, with $A$-action $a(n\ot m)=an\ot m$.\\
If $A$ is commutative and $M,N\in {}_A\Mm^H$, then $M\ot_A N\in {}_A\Mm^H$,
with $H$-coaction
$$\rho_{M\ot_AN}(m\ot n)=m_0\ot n_0\ot m_1n_1.$$

\begin{lemma}\lelabel{2.3}
Let $M\in \Mm^H$, and $N,P\in {}_A\Mm^H$. The isomorphism of vector spaces
$$\phi :\ \Hom^H(N\otimes M, P)\to \Hom^H(M, \HOM(N, P)),~~
\phi(f)(m)(n)=f(n\otimes m),$$
as discussed in \prref{1.2}, induces an isomorphism
$$\phi :\ {}_A\Hom^H(N\otimes M, P)\to \Hom^H(M, {}_A\HOM(N, P)).$$
\end{lemma}

\begin{proof} If $f\in {}_A\Hom^H(N\otimes M, P)$, then $\phi(f)(m)$ is
$A$-linear, for all $m\in M$. If $g\in
\Hom^H(M, {}_A\HOM(N, P))$, then ${{\phi}^{-1}}(g)$ is $A$-linear.
\end{proof}

\begin{corollary}\colabel{2.4} 
Let $I$ be an injective object of ${}_{A}\Mm^H$.
\begin{enumerate}
\item For every $N\in {}_A\Mm^H$,
${}_A\HOM(N, I)$ is an injective object of $\Mm^H$;
\item $I$ is an injective object of $\Mm^H$.
\end{enumerate}
\end{corollary}

\begin{proof} 
(1) follows from \leref{2.3} and the exactness of the functor
$N\otimes (-):\ \Mm^H\to {}_{A}\Mm^H$.\\
(2) By (1), ${}_A\Hom(A, I)$ is an injective object of $\Mm^H$. By \cite[Lemma
1.1]{9}, ${}_A\Hom(A, I)\cong I$ in ${}_{A}\Mm^H$.
\end{proof}

Let $M\in {}_A\Mm^H$. We will use 
the following notation.
\begin{itemize}
\item ${}_A\EXT^q(M,-)$ are the right derived functors of
$${}_A\HOM(M,-):\ {}_A\Mm^H\to \Mm^H;$$
\item ${}_A\Ext^{H^q}(-,-)$ are the right derived functors of
$${}_A\Hom^H:\ {}_A\Mm^H\times {}_A\Mm^H\to \Mm.$$
\end{itemize}
In particular, if $M,N\in {}_A\Mm^H$, then ${}_A\EXT^p(M,N)\in \Mm^H$.

\begin{lemma}\lelabel{2.5} 
For any $N\in {}_{A}\Mm^H$, we have $R^{p}a({\rm co}H,
N)={}_A\Ext^{H^p}(A, N)$.
\end{lemma}

\begin{proof} 
By \leref{2.2}, ${}_A\Hom^H(A, N)={}_A\Hom(A, N)^{{\rm co}H}$. By \cite[Lemma
1.1]{9}, ${}_A\Hom(A, N)^{{\rm co}H}=N^{{\rm co}H}$.
By \coref{2.4} (2), an injective resolution of $N\in{}_{A}\Mm^H$ is also an
injective resolution of $N\in\Mm^H$.
\end{proof}

\begin{proposition}\prlabel{2.6}
Let $M,N\in {}_{A}\Mm^H$. Then we
have a spectral sequence 
\begin{equation}\eqlabel{2.6.1}
R^{p}a({\rm co}H, {}_A\EXT^q(M, N))~~ \Rightarrow~~
{}_A\Ext^{H^{p+q}}(M, N).
\end{equation}
\end{proposition}

\begin{proof}
We have that
$${}_A\Hom^H(M, N)={}_A\HOM(M, N)^{{\rm co}H}.$$
By \coref{2.4}, the functor ${}_A\HOM(M, -)$ takes injective objects of 
${}_A\Mm^H$ to injective objects of $\Mm^H$. 
The result then follows from Grothendieck's spectral sequence for
composite functors.
\end{proof}

\begin{corollary}\colabel{2.7}
Assume that $H$ is cosemisimple, and take $M,N\in {}_{A}\Mm^H$. Then
$${}_A\EXT^q(M, N)^{{\rm co}H}={}_A\Ext^{H^q}(M, N).$$
\end{corollary}

\begin{proof} We know that $\Mm^H$ is a semisimple category. The
result follows from \prref{2.6}.
\end{proof}

\begin{proposition}\prlabel{2.8}
\begin{enumerate}
\item 
For any $I\in {}_{A}\Mm^H$, the functors
${}_A\HOM(I , -)$ and ${}_A\HOM(- , I)$, from  $_{A}\Mm^H$ to  $\Mm^H$,
are left exact.
\item If $I\in {}_{A}\Mm^H$ is injective, then ${}_A\HOM(- , I)$ is
exact.
\end{enumerate}
\end{proposition}

\begin{proof}
 (1) Let $0 \rightarrow M \rTo^{i} N \rTo^{\pi} P
\rightarrow 0$ be an exact
sequence in ${}_{A}\Mm^H$. By \coref{1.7},
$$0 \rightarrow \HOM(P, I) \rightarrow \HOM(N, I) \rightarrow \HOM(M, I)$$
is an exact sequence in $\Mm^H$. It is clear that $f \circ \pi \in
{}_A\HOM(N, I)$ for all
$f \in {}_A\HOM(P, I)$, and $f \circ i \in {}_A\HOM(M, I)$ for all $f \in
{}_A\HOM(N, I)$. It follows that
$$0 \rightarrow {}_A\HOM(P, I) \rightarrow {}_A\HOM(N, I) \rightarrow {}_A\HOM(M, I)$$
is an exact sequence in $\Mm^H$. So the functor ${}_A\HOM(- , I)$ is left
exact. In a similar way, we show that
the functor ${}_A\HOM(I , -)$ is left exact.

(2) Let $f \in {}_A\HOM(M, I)$ and let $V$ be a finite-dimensional
$H$-subcomodule of $\HOM(M, I)$
containing $f$. The $k$-linear map
$$i \otimes V: M\otimes V \rightarrow N\otimes V$$
is a monomorphism of relative $(A, H)$-Hopf modules. The map
$$F:\ M\ot V\to I,~~F(m\ot v)=v(m)$$
is $A$-linear. 
As in the proof of $(1)
\Rightarrow (2)$ in \prref{1.5}, we can show that  $F$ is 
$H$-colinear, hence $F$ is a morphism in ${}_A\Mm^H$.
Since $I$ is injective, there exists a morphism
$G:\ N\otimes V \rightarrow
I$ in ${}_A\Mm^H$ such that
$G \circ (i \otimes V)=F$. The map
$$g:\ N\to I,~~g(n)=G(n\ot f)$$
is right $A$-linear, and it follows from \prref{1.5} that
$g$ is rational, hence $g\in {}_A\HOM(N,I)$. Finally
$$f(m)=F(m \otimes f)=G \circ (i \otimes V)(m \otimes f)= G(i(m) \otimes
f)=(g \circ i)(m),$$
and it follows that ${}_A\HOM(N, I) \rightarrow {}_A\HOM(M, I)$ is
surjective.
\end{proof}

\begin{proposition}\prlabel{2.9}
Let $A$ be left noetherian and $M\in {}_A\Mm^H$ finitely generated as a left $A$-module.
If $I\in {}_A\Mm^H$ is injective, then ${}_A\Ext^p(M, I)=0$, for all $p>0$.
\end{proposition}

\begin{proof}
Since $M$ is finitely generated, there exist a finite dimensional
$H$-comodule
$V_0$ and an epimorphism $p_0:\ P_0=A\otimes V_0
\rightarrow M$ in ${}_A\Mm^H$, by \cite[Prop. 4.1]{9}. 
Then $K=\Ker(p_0)$ is a subobject of $P_0$ in ${}_A\Mm^H$. $P_0\in {}_A\Mm$
is finitely generated, and $A$ is left noetherian, so
$K\in {}_A\Mm$ is also finitely generated. So we can find a finite dimensional
$H$-comodule $V_1$ and an epimorphism $p_1:\ P_1=A\ot V_1\to K$, and
we have that $\im(p_1)=K=\Ker(p_0)$. Repeating this construction, we find
an $A$-free resolution $P_\star$ of $M$ in ${}_A\Mm^H$,
$$\cdots \rightarrow P_i= A\otimes V_i \rightarrow ...
\rightarrow  P_1=
A\otimes V_1 \rightarrow  P_0= A\otimes V_0 \rightarrow M \rightarrow 0,$$
with each $V_i$ a finite dimensional $H$-comodule. For each $p>0$, we have
$${}_A\Ext^p(M,I)=H^p({}_A\Hom(P_{*}, I)).$$
From \leref{2.2}, we know that
$${}_A\HOM(M, I)={}_A\Hom(M, I)~~{\rm and}~~
{}_A\HOM(P_i, I)={}_A\Hom(P_i, I),$$
for all $i\geq 0$.
On the other hand, $P_{\star}$ is an acyclic complex in ${}_{A}\Mm^H$.
We deduce from \coref{2.4} and \prref{2.8} that ${}_A\HOM(P_{\star}, I)$ is an injective
resolution
of ${}_A\HOM(M, I)$ in $\Mm^H$, and it follows that $H^p({}_A\HOM(P_{\star},
I))=0$ for all $p>0$.
\end{proof}

\begin{corollary}\colabel{2.10}
Let $A$ be left noetherian. Take $M,N\in {}_A\Mm^H$, with $M$ finitely
generated as an $A$-module and $E^*=\{E^i\}$ an injective resolution of $N$ in
${}_A\Mm^H$. Then for all $p\geq 0$
\begin{equation}\eqlabel{2.10.1}
{}_A\Ext^p(M, N)={}_A\EXT^p(M, N)=H^p({}_A\Hom(M, E^*)).
\end{equation}
\end{corollary}

\begin{proof} 
For all $p\geq 0$, we have that 
$${}_A\Ext^p(M, N)=H^p({}_A\Hom(P_\star,N))=H^p({}_A\HOM(P_{\star}, N)),$$
where $P_\star$ is the $A$-free resolution of $M$ constructed in \prref{2.9}.
${}_A\HOM(P_{\star}, N)$ is a complex in $\Mm^H$ which induces on each
$H^p({}_A\HOM(P_{\star}, N))$ $  ={}_A\Ext^p(M,
N)$ a structure of $H$-comodule, so ${}_A\Ext^p(M, -)$ is a cohomological
functor from ${}_{A}\Mm^H$ to
$\Mm^H$ and, by \prref{2.9}, ${}_A\Ext^p(M, I)=0$ for all $p>0$ 
if $I\in {}_A\Mm^H$ is injective. Clearly the same property holds for
${}_A\EXT^p(M, -)$. By \leref{2.2}, we have that
$${}_A\Ext^0(M, N)={}_A\Hom(M, N)={}_A\HOM(M, N)=R^0({}_A\HOM(M, -))(N)$$ 
in ${}_A\Mm^H$. It follows that
${}_A\Ext^p(M, -)$ and ${}_A\EXT^p(M, -)$ coincide on ${}_A\Mm^H$ for $p\geq
0$, and we obtain the first
equality of \equref{2.10.1}. The second one follows after we observe that ${}_A\EXT^p(M, N)=H^p({}_A\HOM(M, E^*)$
and, by \leref{2.2}, ${}_A\HOM(M,
E^p)={}_A\Hom(M, E^p)$, for all $p \geq 0$.
\end{proof}

\begin{proposition}\prlabel{2.11} 
Let $A$ be left noetherian. Take $M,N\in {}_A\Mm^H$,
with $M$ finitely generated as a left $A$-module. Then
we have a spectral sequence
$$R^{p}a({\rm co}H, {}_A\Ext^q(M, N))~~\Rightarrow~~ {{}_A\Ext^H}^{p+q}(M, N).$$
\end{proposition}

\begin{proof}
By \leref{2.2}, we know that ${}_A\Hom(M, N)^{{\rm co}H}= {}_A\Hom^H(M, N)$.
By \prref{2.8}, the
functor ${}_A\Hom(M , -)$ takes injective objects of ${}_{A}\Mm^H$ to 
injective objects of $\Mm^H$. 
Now $R^q({}_A\Hom(M , -))(N)={}_A\Ext^q(M , N)$ for every $q\geq 0$,
by \coref{2.10}.
The result then follows from the
Grothendieck spectral sequence for composite functors.
\end{proof}

\begin{corollary}\colabel{2.12}
Assume that $H$ is cosemisimple, and
that $A$ is left noetherian. Take $M,N\in {}_A\Mm^H$,
with $M$ finitely generated as a left $A$-module. Then
$${}_A\Ext^q(M, N)^{{\rm co}H}={{}_A\Ext^H}^q(M, N).$$
\end{corollary}

\begin{proof} We know that $\Mm^H$ is a semisimple category, so the
result follows from \coref{2.10}.
\end{proof}

With notation and assumptions as in \coref{2.12}, it follows that if
$M\in {}_{A}\Mm^H$ is finitely
generated and projective in ${}_{A}\Mm$, then $M$ is also projective in
${}_{A}\Mm^H$.

\begin{lemma}\lelabel{2.13}
Let $A$ and $H$ be commutative. Let take $M, N\in {}_A\Mm^H$.
Then ${}_A\HOM(M, N)\in {}_A\Mm^H$. A fortiori ${}_A\EXT^p(M , N)\in {}_A\Mm^H$.
\end{lemma}

\begin{proof} 
By \leref{2.2}, ${}_A\HOM(M, N)$ is a a right $H$-comodule. For $a\in A$,
we consider the $k$-linear map
$$L(a):\ M\rightarrow M,~~L(a)(m)=am.$$ 
Then for all $m\in M$, we have that
\begin{eqnarray*}
&&\hspace*{-2cm}
(L(a)(m_0))_0 \otimes (L(a)(m_0))_1S(m_1)=(am_0)_0 \otimes (am_0)_1S(m_1)\\
&=& a_0m_0 \otimes a_1m_1S(m_2)=a_0m \otimes a_1=L(a_0)(m) \otimes a_1,
\end{eqnarray*}
so $L(a)_0 \otimes L(a)_1 = L(a_0) \otimes a_1$, and $L\in {}_A\HOM(M, M)$.
For $f\in {}_A\HOM(M,N)$, we now set $af=f\circ L(a)$. It follows from
\prref{1.5} that $af\in \HOM(M,N)$, and it is clear that $af$ is left $A$-linear.
Hence ${}_A\HOM(M,N)$ is a left $A$-module. Let us finally check the compatibility
relation between the action and coaction on ${}_A\HOM(M,N)$. For all
 $f\in {}_A\HOM(M, N)$, $m \in M$ and $a \in A$, we have
\begin{eqnarray*}
&&\hspace*{-2cm} 
((af)_0 \otimes (af)_1)(m)=((af)(m_0))_0 \otimes ((af)(m_0))_1S(m_1)\\
&=&a_0(f(m_0)_0) \otimes a_1(f(m_0)_1)S(m_1)\\
&=&a_0(f(m_0)_0) \otimes
a_1(f(m_0)_1S(m_1))\\
&=& a_0(f_0(m)) \otimes a_1f_1\\
&=&(a_0f_0)(m) \otimes a_1f_1=(a_0f_0 \otimes
a_1f_1)(m).
\end{eqnarray*}
\end{proof}

Let $A$ be commutative, and take
$M,N\in {}_{A}\Mm^H$. 
By \cite[Lemma 1.1]{9}, $M\otimes_AN\in {}_{A}\Mm^H$. The action and coaction
are given by the formulas
$$a(m\otimes n)=am \otimes n=m\otimes an;$$
$${\rho}_{M\otimes_AN}(m\otimes n)=m_0 \otimes n_0\otimes m_1n_1.$$

\begin{proposition}\prlabel{2.14}
Let $A$ and $H$ be commutative, take $M,N,P\in {}_{A}\Mm^H$, and
consider the natural $k$-isomorphism
$$\phi:\ {}_A\Hom(M\otimes_AN, P)\to {}_A\Hom(M , {}_A\Hom(N, P)),~~
\phi(f)(m)(n)=f(m\otimes n).$$
\begin{enumerate}
\item If $f\in {}_A\Hom(M\otimes_A N, P)$ is $H$-colinear, then $\phi(f)(m)\in
{}_A\HOM(N, P)$,
for every $m \in M$; furthermore $\phi(f)$ is $H$-colinear;
\item $\phi$ induces a $k$-isomorphism
$$\phi:\ {}_A\Hom^H(M\otimes_A N, P)\to
{}_A\Hom^H(M, {}_A\HOM(N, P));$$
\item If $N$ is flat as a left $A$-module, then ${}_A\HOM(N, -)$ preserves the
injective objects of ${}_{A}\Mm^H$.
\end{enumerate}
\end{proposition}

\begin{proof} (1) and (2): an easy adaptation of the proof of (1) and (2) in \prref{1.2}.

(3) If $I\in {}_{A}\Mm^H$ is injective, then the functor ${}_A\Hom^H(-
, I)$ is exact. $N$ is flat as a left $A$-module, so
$-\otimes_AN$ is an exact endofunctor of
${}_{A}\Mm^H$. It then follows from (2)
that the functor ${}_A\Hom^H(- , {}_A\HOM(N , I))$ is exact.
\end{proof}

\begin{proposition}\prlabel{2.15} 
Let $A$ and $H$ be commutative, and take $M,N,P\in {}_{A}\Mm^H$.
If $N$ is flat as a left $A$-module, then we have a spectral sequence
$${{}_A\Ext^H}^p(M , {}_A\EXT^q(N, P))~~\Rightarrow~~ {{}_A\Ext^H}^{p+q}(M\otimes_AN, P).$$
\end{proposition}

\begin{proof} The functors ${}_A\Hom(M\otimes_AN, -)$ and
${}_A\Hom(M, {}_A\HOM(N, -))$ coincide on ${}_{A}\Mm^H$, by \prref{2.14} (2).
${}_A\HOM(M , -)$ preserves the injectives of ${}_{A}\Mm^H$, by \prref{2.14} (3).
\end{proof}

\section{The functor ${}_B\HOM(A , -)$}\selabel{3}
Recall that $\varsigma \in H^*$ is called a left
integral on $H$ if
$h^{*}\phi=h^*(1)\phi$ for all
$h^*\in H^*$.
Throughout this Section, we assume that $H$ is cosemisimple, which is equivalent to the existence of a left integral
$\phi$ on $H^*$ such that $\phi(1)=1$ (see e.g. \cite{27}). For every $M\in \Mm^H$,
we then have an $H$-colinear epimorphism (see \cite[Prop. 1.5]{2})
$$p_M:\ M\to M^{{\rm co}H},~~p_M(m)=\phi(m_1)m_0.$$

$M\in \Mm^H$ is called ergodic if $M^{{\rm co}H}=0$. A subcomodule of an ergodic comodule
is ergodic, and, for every $M\in \Mm^H$, 
$M/M^{{\rm co}H}$ is ergodic. Let $M_{{\rm co}H}$ be the maximal ergodic subcomodule
of $M$. It is obvious that $M^{{\rm co}H}\cap M_{{\rm co}H}=0$, and we have

\begin{lemma}\lelabel{3.1} 
Let $H$ be a cosemisimple Hopf algebra. Then for all $M\in \Mm^H$,
$$M=M^{{\rm co}H}\oplus M_{{\rm co}H}$$
as $H$-comodules.
\end{lemma}

The decomposition of \leref{3.1} is functorial in the following
sense. If $f:\ M\to M'$ is $H$-colinear, then
$f(M^{{\rm co}H})\subseteq M'^{{\rm co}H}$ and $f(M_{{\rm co}H})\subseteq M'_{{\rm co}H}$. In
particular, the projection $p_M:\ M\to M^{{\rm co}H}$ is $H$-colinear, and
$f\circ p_M=p_{M'}\circ f$.

Let $M$ be a $B$-module, and let $H$ coact trivially on $M$. In particular,
$H$ coacts trivially on $B$, $B$ is an $H$-comodule algebra, and
$M$ is a relative $(B,H)$-Hopf module.

Take $M\in {}_A\Mm^H$. For $b\in B=A^{{\rm co}H}$, the map $f_b\in \End(M)$
given by 
$f_b(m)= bm$ is $H$-colinear, so $f_b\circ p_M=p_M\circ
f_b$. It follows that $f_b(M^{{\rm co}H})\subseteq M^{{\rm co}H}$ and
$f_b(M_{{\rm co}H})\subseteq
M_{{\rm co}H}$, that is, $M^{{\rm co}H}$ and $M_{{\rm co}H}$ are $B$-submodules (hence
$(B, H)$-Hopf submodules) of $M$ and $p_M$ is $B$-linear (hence a morphism
of $(B, H)$-Hopf modules).

Recall from \cite[Lemmas 2.1 and 2.2]{2} that ${}_B\HOM(A,M)\in {}_A\Mm^H$.
The left $A$-action is given by the formula
$$(af)(a')=f(a'a).$$

For every  $M$ (resp. $N$) in ${}_{B}\Mm$ (resp. in ${}_{A}\Mm^H$),
${}_BE(M)$ (resp. ${{}_AE^H}(N)$) will be the injective hull of $M$
in ${}_{B}\Mm$ (resp. of $N$ in
${}_A\Mm^H$).

\begin{lemma}\lelabel{3.2}
Let $A$ be an $H$-comodule algebra.
\begin{enumerate} 
\item If $M\in {}_{B}\Mm$, then $A\otimes_BM\in {}_{A}\Mm^H$.
\item
\begin{enumerate}
\item For $M\in {}_{B}\Mm$ and $N\in {}_{A}\Mm^H$, we have an
isomorphism of $k$-vector spaces
$${}_A\Hom^H(A\otimes_BM, N)\cong {}_B\Hom(M , N^{{\rm co}H});$$
\item For $M\in {}_{B}\Mm$ and $N\in {}_{B}\Mm^H$, we have an
isomorphism of $k$-vector spaces
$${}_B\Hom^H(M, N)\cong {}_B\Hom(M , N^{{\rm co}H}).$$
\end{enumerate}
\end{enumerate}
\end{lemma}

\begin{proof} (1) is obvious. (2a) follows from the fact that we have
a pair of adjoint functors $(A\ot_B -,(-)^{{\rm co}H})$ between
${}_{A}\Mm^H$ and ${}_{B}\Mm$. (2b) follows after we take $A=B$ in (2a),
with trivial coaction on $A$.
\end{proof}

Also recall the following results from \cite[Theorem 2.3, Cor. 2.4 and 2.5]{2}.

\begin{proposition}\prlabel{3.3}
Let $A$ be an $H$-comodule algebra, and assume that $H$ is cosemisimple.
Take $N\in{}_B\Mm$ and $M\in {}_A\Mm^H$.
\begin{enumerate}
\item The map
$$\phi:\ {}_A\Hom^H(M ,{}_B\HOM(A , N))\rightarrow {}_B\Hom(M^{{\rm co}H} , N),~~\phi(f)(p_M(m))=f(m)(1)$$ is an
isomorphism of $k$-vector spaces;
\item the map 
$$F: {}_B\HOM(A , N)^{{\rm co}H}\rightarrow N,~~F(f)=f(1)$$
is an isomorphism of $B$-modules;
\item if $I\in {}_B\Mm$ is injective, then ${}_B\HOM(A,I)\in {}_A\Mm^H$
is injective.
\end{enumerate}
\end{proposition}

\begin{theorem}\thlabel{3.6}
Let $A$ be an $H$-comodule algebra, and assume that $H$ is cosemisimple.
\begin{enumerate}
\item If $N\in {}_B\Mm$ and $M$ is an $(A , H)$-Hopf
submodule of ${}_B\HOM(A , N)$, then
$M^{{\rm co}H}=0$ implies $M=0$;
\item if $M\rightarrow N$ is an essential monomorphism in ${}_B\Mm$,
then ${}_B\HOM(A , M)\rightarrow {}_B\HOM(A , N)$ is an essential monomorphism
in ${}_A\Mm^H$;
\item if $N\in {}_{B}\Mm$, then ${{}_AE^H}({}_B\HOM(A , N))\cong{}_B\HOM(A , {}_BE(N))$;
\item if $N\in {}_{B}\Mm$, then $({{}_AE^H}({}_B\HOM(A , N)))^{{\rm co}H}\cong
{}_BE(N)$.
\end{enumerate}
\end{theorem}

\begin{proof} For any subset $T$ of ${}_B\HOM(A , M)$, set
$T(1)=\{f(1)~|~f\in T\}$.

(1) If $M^{{\rm co}H}=0$, then ${}_A\Hom(M^{{\rm co}H} , N)=0$,
and, by \prref{3.3}(1), ${}_A\Hom^H(M ,{}_B\HOM(A , N))=0$. Hence
the inclusion map $M\to {}_B\HOM(A , N)$ is the zero map, hence $M=0$.

(2) If $L$ is a nonzero $(A , H)$-Hopf submodule of $_B\HOM(A, N)$, then by
(1), $L^{{\rm co}H}$ is a nonzero
$B$-submodule of ${}_B\HOM(A, N)$. By \prref{3.3}(2), this means that $L(1)$ is a
nonzero $B$-submodule of $N$, 
so $L(1)\cap M\not=0$. But $L(1)\cap M=({}_B\HOM(A , N)\cap L)(1)$; so $L$
meets ${}_B\HOM(A , N)$ nontrivially.

(3) By (2), ${}_B\HOM(A , N)\rightarrow {}_B\HOM(A ,{}_BE(N))$ is an
essential monomorphism in ${}_{A}\Mm^H$. 
But, by (2), ${}_B\HOM(A ,{}_BE(N))$ is an injective object of ${}_{A}\Mm^H$.

(4) It follows from (3) that
${{}_AE^H}({}_B\HOM(A , N))={}_B\HOM(A,{}_BE(N))$, and from \prref{3.3}(3)
that
${}_B\HOM(A ,{}_BE(N))^{{\rm co}H}={}_BE(N)$.
\end{proof}

By (1), the nonzero subobjects of ${}_B\HOM(A,N)$ in ${}_A\Mm^H$ contain
nonzero coinvariants. We will see below that this - rather strong -
property implies that $M$ is an
essential extension of $AM^{{\rm co}H}$.\\

Let $H^*$ be the linear dual of $H$, and consider the smash product $A\# H^*$
(see e.g. \cite{12}). Then we have a functor ${}_A\Mm^H\to {}_{A\# H^*}\Mm$,
and, conversely, a left $A\#H^*$-module which is rational as an $H^*$-module can be regarded as a relative $(A,H)$-Hopf module.

\begin{corollary}\colabel{3.7}
Take $N\in {}_{B}\Mm$ and let $M\neq 0$ be a subobject of ${}_B\HOM(A , N)$
in ${}_A\Mm^H$. Take $m\in M$.
\begin{enumerate}
\item $M$ is an essential extension of $AM^{{\rm co}H}$ in ${}_A\Mm^H$.
\item If $p_M(am)=0$ for all $a\in A$, then $m=0$.
\end{enumerate}
\end{corollary}

\begin{proof} 
(1) Let $L\neq 0$ be a subobject of $M$ in ${}_A\Mm^H$. By
\thref{3.6}(1), $L^{{\rm co}H}\not=0$
and $L^{{\rm co}H}\subseteq M^{{\rm co}H}$, so $L\cap AM^{{\rm co}H}\not=0$.

(2) By \cite[p. 247]{11}, $A\#H^*$ is isomorphic as a left $H^*$-module to
$H^* \otimes A$. So each element
of $A\#H^*$ can be written as a finite sum $\sum_i {h_i^*}a_i$, with
$h_i^* \in H^*$ and $a_i \in
A$. Consider the $A\#H^*$-submodule $P$ of $M$ generated by $m$. If $m\neq 0$
then, by
\thref{3.6}(1), $P$ contains a nonzero coinvariant element $y=\sum
h_i^*a_im$. But $p_M(y)=y$ while $p_M(\sum h_i^*a_im)=\sum
h_i^*p_M(a_im)=0$, since $p_M$ is
$H^*$-linear. So $y=0$, which is a contradiction. We conclude that $m=0$.
\end{proof}

For $M\in {}_{A}\Mm^H$, we set 
$${}^{\bullet}M=\{m\in M ~|~ p_M(am)=0\quad \hbox{for all}\quad a\in A\}.$$
Note that if $M\in {}_A\Mm^H$ is simple, and $M^{{\rm co}H}\not=0$,
then ${}^{\bullet}M=0$. Indeed,
if ${}^{\bullet}M=M$, then $p_M(m)=0$ for every $m\in M$, so
$M=M_{{\rm co}H}$, hence $M^{{\rm co}H}=0$ which is a
contradiction.

\begin{lemma}\lelabel{3.8}
 Let $M\in {}_{A}\Mm^H$ and consider the natural transformation
 (see \cite[Prop. 2.7]{2})
$$\nu_M:\ M\rightarrow {}_B\HOM(A , M^{{\rm co}H}),~~
\nu_M(m)(a)=p_M(am).$$
\begin{enumerate}
\item ${}^{\bullet}M=\Ker \nu_M$;
\item if $f : M \rightarrow M'$ is a morphism in ${}_A\Mm^H$ 
then $f({}^{\bullet}M) \subset {}^{\bullet}M'$;
\item  ${}^{\bullet}(M/{}^{\bullet}M)=0$;
\item if $M$ is a subobject of $N$ in ${}_A\Mm^H$, then ${}^{\bullet}N\cap M={}^{\bullet}M$;
\item if ${}^{\bullet}M=0$, then $\nu_M$ is an essential monomorphism in $_{A}\Mm^H$.
\end{enumerate}
\end{lemma}

\begin{proof}
(1) is obvious, and (2) follows from the fact that $f\circ p_M=
p_{M'}\circ f$.\\
(3) Observe that $({}^{\bullet}M)^{{\rm co}H}=0$. Hence $(M/{}^{\bullet}M)^{{\rm co}H}=M^{{\rm co}H}$. As
$\Ker({\nu}_M)={}^{\bullet}M$, the map $\nu_M$ factorizes through
$$\ol{\nu}_M:\ M/{}^{\bullet}M\to {}_B\HOM(A , M^{{\rm co}H}).$$
Now ${}_B\HOM(A , M^{{\rm co}H})\cong {}_B\HOM(A,(M/{}^\bullet M)^{{\rm co}H})$,
so it follows that $\ol{\nu}_M=\nu_{M/{}^\bullet M}$.
Take $m\in M$ such that the corresponding $[m]\in M/{}^\bullet M$
is in ${}^{\bullet}(M/{}^{\bullet}M)$. It follows from (1) that
$\nu_{M/{}^{\bullet}M}([m])=0$, hence $m\in {}^{\bullet}M$, and
it follows that ${}^{\bullet}(M/{}^{\bullet}M)=0$.\\
(4) follows from the definition of ${}^{\bullet}M$, and the fact that
the restriction of $p_N$ to $M$ is $p_M$.\\
(5) Assume that ${}^{\bullet}M=0$, and identify $M$ with $\nu_M(M)$.
If $L\neq 0$ is a subobject of ${}_B\HOM(A,M^{{\rm co}H})$, then, by
\thref{3.6}, $L^{{\rm co}H}\neq 0$, and, by \prref{3.3}(2), 
$M^{{\rm co}H}={}_B\HOM(A,M^{{\rm co}H})^{{\rm co}H}$, so
$L\cap M\neq 0$.
\end{proof}

\leref{3.8} can be used to characterise injective objects in ${}_A\Mm^H$
of the form ${}_B\HOM(A,I)$, with $I\in {}_B\Mm$ injective.

\begin{theorem}\thlabel{3.9}
\begin{enumerate} 
\item If $E\in {}_{A}\Mm^H$ is injective and
${}^{\bullet}E=0$, then
$E^{{\rm co}H}\in {}_B\Mm$ is injective and 
$E\cong {}_B\HOM(A , E^{coH})$ in ${}_A\Mm^H$.
\item If $M\in {}_{A}\Mm^H$ with ${}^{\bullet}M=0$, then 
${}_AE^H(M)\cong {}_B\HOM(A, {}_BE(M^{{\rm co}H}))$  in ${}_A\Mm^H$.
\end{enumerate}
\end{theorem}

\begin{proof} 
(1) Set $E'={}_BE(E^{{\rm co}H})$. Then $E^{{\rm co}H} \to E'$ is
an essential monomorphism of
$B$-modules, so, by \thref{3.6}(2), ${}_B\HOM(A , E^{{\rm co}H}) \to {}_B\HOM(A ,
E')$ is an essential monomorphism in ${}_A\Mm^H$.
Since $E$ is injective in ${}_{A}\Mm^H$, we have that
$$E\cong {}_B\HOM(A , E^{{\rm co}H}) \cong {}_B\HOM(A , E').$$
By \prref{3.3}(2), ${}_B\HOM(A , E')^{{\rm co}H} \cong
E'$ is an injective $B$-module, so $E^{{\rm co}H}\cong E'$ is an injective
$B$-module.

(2) Set $E={}_BE(M^{{\rm co}H})$. By \prref{3.3}(3), ${}_B\HOM(A , E)\in {}_{A}\Mm^H$ is injective,
and by \leref{3.8}(5) and
\thref{3.6}(2), we have essential monomorphisms
$$M \rightarrow {}_B\HOM(A , M^{{\rm co}H}) \rightarrow {}_B\HOM(A , E)$$
in ${}_A\Mm^H$, and it
follows that ${}_AE^H(M)\cong{}_B\HOM(A , E)$ in ${}_A\Mm^H$.
\end{proof}

\begin{remark}\relabel{3.9bis}
It follows from \coref{3.7}(2) that we have the following converse of
\thref{3.9}(1):
if an injective object $E\in {}_A\Mm^H$ is isomorphic to 
${}_B\HOM(A , E^{coH})$ in ${}_A\Mm^H$, then ${}^\bullet E=0$.
\end{remark}

It is well-known that ${}_A\Mm$ has an injective cogenerator $I$,
and it follows from \cite[Prop. 1, Theorem 3]{20} that $I\ot H$
is an injective cogenerator of ${}_A\Mm^H$.\\
If ${}_A\Mm^H$ has an injective cogenerator $C$ with ${}^\bullet C=0$,
then it follows that ${}^\bullet M=0$, for every $M\in {}_A\Mm^H$,
by \leref{3.8}(4). In this case, we will say that ${}_A\Mm^H$ satisfies
the condition $(\alpha)$.

\begin{proposition}\prlabel{3.10}
Assume that ${}_{A}\Mm^H$ satisfies condition $(\alpha)$.
\begin{enumerate}
\item Every injective object of ${}_{A}\Mm^H$ is isomorphic to
${}_B\HOM(A , I)$, for some injective left $A$-module $I$;
\item For $M\in {}_{B}\Mm$ and $N\in {}_{A}\Mm^H$, we have
$${{}_A\Ext^H}^p(A\otimes_BM , N)\cong {}_B\Ext^p(M , N^{{\rm co}H}),$$
for all $p\geq 0$.
\end{enumerate}
\end{proposition}

\begin{proof} 
(1) follows immediately from \thref{3.9} and \reref{3.9bis}.

(2) Let $E^*=\{E^i\}$ be an injective resolution of $N$ in $_{A}\Mm^H$. By \leref{3.2}(2),
$${}_A\Hom^H(A\otimes_BM , E^i) \cong {}_B\Hom(M , {E^i}^{{\rm co}H}),$$
for every $i$, so we have that
\begin{equation}\eqlabel{3.10.1}
{{}_A\Ext^H}^p(A\otimes_BM , N)=H^p({}_B\Hom(M , {E^*}^{{\rm co}H})).
\end{equation}
for all $ p\geq 0$.
But the functor $(-)^{{\rm co}H}$ is exact and by \thref{3.9}(1), each
${E^i}^{{\rm co}H}$ is injective in ${}_{B}\Mm$. Hence $\{{E^i}^{{\rm co}H}\}$ is an
injective resolution of
$N^{{\rm co}H}$ in ${}_B\Mm$, and the right hand side of \equref{3.10.1}
is ${}_B\Ext^p(M , N^{{\rm co}H})$.
\end{proof}

\begin{lemma}\lelabel{3.11}
Let $A$ be noetherian, and $\{E^i~|~ i\in I\}$ be a set of
injective
$B$-modules. We have the following isomorphism in  ${}_A\Mm^H$:
$$E=\bigoplus_{i\in I} {}_B\HOM(A , E^i)\cong {}_B\HOM\Bigl(A, \bigoplus_{i\in I} E_i
\Bigr).$$
\end{lemma}

\begin{proof} 
$E\in {}_A\Mm^H$ is injective, by \prref{3.3}(3), and
${}^{\bullet}_B\HOM(A,E^i)=0$ for all $i\in I$, by \coref{3.7}(2), hence
${}^{\bullet}E=0$. $E^{{\rm co}H}= \bigoplus_{i\in I} E^i$, by
\thref{3.6}. We have seen in \leref{3.8}(5) that $\nu_E$ is an essential
monomorphism in ${}_A\Mm^H$ and, since $E\in {}_A\Mm^H$ is injective,
$\nu_E$ is an isomorphism in ${}_A\Mm^H$.
\end{proof}

\begin{lemma}\lelabel{3.12} 
Let $I\in {}_B\Mm$ be injective, and take $M\in {}_{A}\Mm^H$.
Assume that ${}^{\bullet}M=0$ and that
$f:\ M\to {}_B\HOM(A, I)$ is an essential monomorphism in ${}_{A}\Mm^H$. Then
$${M^{{\rm co}H}}\to {}_BHOM(A , I)^{{\rm co}H}=I$$ 
is an essential monomorphism in ${}_B\Mm$.
\end{lemma}

\begin{proof} 
By \leref{3.8}(5), $\nu_M$ is an essential monomorphism in
${}_{A}\Mm^H$ and, by \prref{3.3}(3), ${}_B\HOM(A ,
I)$ is injective in ${}_{A}\Mm^H$. So there exists a morphism
$h :\ {}_B\HOM(A ,M^{{\rm co}H}) \to {}_B\HOM(A , I)$ in ${}_{A}\Mm^H$
such that $f=h \circ \nu_M$. Let $L$ be
a $B$-submodule of $I$ such that
$L\cap M^{{\rm co}H}=0$. Then
$${}_B\HOM(A , M^{{\rm co}H}) \cap {}_B\HOM(A , L)=0$$
and ${}_B\HOM(A , L)$ is a relative $(A , H)$-Hopf submodule of ${}_B\HOM(A , I)$. If
${}_B\HOM(A , L)\neq 0$ then ${}_B\HOM(A , L)$
meets $M=\mu_M(M)$ nontrivially because $\nu_M(M)\subseteq {}_B\HOM(A ,
M^{{\rm co}H})$. This is impossible, so
${}_B\HOM(A , L)=0$. We deduce from \prref{3.3}(2) that $0={}_B\HOM(A , L)^{{\rm co}H}=L$.
\end{proof}

Now we are ready to show that the functor $(-)^{{\rm co}H}$ takes minimal
injective resolutions of ${}_{A}\Mm^H$ to minimal injective
resolutions of ${}_{B}\Mm$.

\begin{proposition}\prlabel{3.13}
Assume that ${}_A\Mm^H$ satisfies condition $(\alpha)$, and that
$A$ and $B$ are noetherian. Take $M\in {}_{A}\Mm^H$, and let
$\{{{_AE^H}}^i(M)\}$ be the minimal injective resolution of $M$ in
${}_{A}\Mm^H$ and $\{{}_BE^i(M^{{\rm co}H})\}$ the minimal injective
resolution of $M^{{\rm co}H}$ in ${}_{B}\Mm$. Then
$({_AE^H}^i(M))^{{\rm co}H}={}_BE^i(M^{{\rm co}H})$, for all $i$.
\end{proposition}

\begin{proof} 
Set $E^i={{}_AE^H}^i(M)$ and $K^i=\Ker(E^i
\to E^{i+1})$, for all $i\geq 0$. It follows from \thref{3.9}(1) that
$I^i={E^i}^{{\rm co}H}$ is an injective $B$-module and $E^i={}_B\HOM(A , I^i)$.
Since $H$ is cosemisimple, the
sequence
$$I^0 \rightarrow I^1 \rightarrow \cdots \rightarrow I^i \rightarrow \cdots$$
is exact in ${}_B\Mm$ and $(K^i)^{{\rm co}H}=\Ker(I^i \rightarrow I^{i+1})$.
Since ${}^{\bullet}(K^i)=0$ and
$K^i \to E^i={}_B\HOM(A , I^i)$ is an essential monomorphism in ${}_A\Mm^H$,
${K^i}^{{\rm co}H} \to I^i={}_B\HOM(A , I^i)^{{\rm co}H}$ is an essential
monomorphism of $B$-modules, by \leref{3.12}, so
$\{I^i\}$ is a minimal injective resolution of ${K^0}^{{\rm co}H}=M^{{\rm co}H}$ in
${}_B\Mm$.
\end{proof}

\begin{theorem}\thlabel{3.14}
Assume that ${}_A\Mm^H$ satisfies condition $(\alpha)$, and that $A$ and $B$ are
noetherian, with $B$ commutative. Take $M\in {}_{A}\Mm^H$. For every
$P\in \Spec(B)$, let $\mu_i(P,
M^{{\rm co}H})$ be the number of times that ${}_BE(B/P)$ occurs in 
${}_BE^i(M^{{\rm co}H})$. Then
$${{}_AE^H}^i(M)=\bigoplus_{P\in \Spec(B)}{{}_AE^H}(A/PA)^{\mu_i(P ,
M^{{\rm co}H})}.$$
\end{theorem}

\begin{proof} 
By \thref{3.9}(1) and \prref{3.13}, $({{}_AE^H}^i(M))^{{\rm co}H}$ is
$B$-injective and ${{}_AE^H}^i(M)\cong {}_B\HOM(A , ({{_AE^H}}^i(M))^{{\rm co}H})=
{}_B\HOM(A ,
{}_BE^i(M^{{\rm co}H}))$ in ${}_A\Mm^H$. So by the definition
of $\mu_i$ and \leref{3.11},
${{}_AE^H}^i(M)$ is the direct sum over $P\in spec(B)$ of $\mu_i(P ,
M^{{\rm co}H})$ copies of ${}_B\HOM(A ,
_BE(B/P))$. But $(A/PA)^{{\rm co}H}=B/P$, so, by \thref{3.9}(2), 
${}_B\HOM(A , _BE(B/P))\cong {}_AE^H(A/PA)$ in ${}_A\Mm^H$.
\end{proof}

\begin{lemma}\lelabel{3.15}
Take $M\in {}_{B}\Mm$, $N\in {}_{A}\Mm^H$ and $V\in \Mm^H$ finite dimensional.
\begin{enumerate}
\item Assume that H has the symmetry property.
We have the following isomorphisms in $\Mm^H$:
$${}_A\Hom(A\otimes V , N)\cong \Hom(V, N)\cong V^*\otimes N\cong N\otimes
V^*.$$
Consequently, $(N\otimes V^*)^{{\rm co}H}$ and $\Hom^H(V , N)$
are isomorphic as vector spaces.
\item If $B$ be commutative, then the map 
$$\phi :\ {}_B\Hom^H(V\otimes A ,
M)\to \Hom^H(V, {}_B\HOM(A , M)),~~\phi(f)(v)(a)=f(v\otimes a)$$ 
is an isomorphism of $k$-vectorspaces.
\end{enumerate}
\end{lemma}

\begin{proof}
(1) The first two isomorphisms are well-known; the third one is a
consequence of the symmetry property.\\
(2) is a direct consequence of \prref{1.2}(3).
\end{proof}

\begin{remark}\relabel{3.15bis}
If $A$ and $H$ are commutative, then the isomorphisms in \leref{3.15}(1)
are left $A$-linear, and therefore $(N\otimes V^*)^{{\rm co}H}$ and $\Hom^H(V , N)$
are isomorphic left $B$-modules.
\end{remark}

The condition that ${}_A\Mm^H$ satisfies the condition $(\alpha)$ is quite
restrictive; it implies that the coinvariants functor
${}_A\Mm^H\to {}_B\Mm$ preserves injectivity (this follows from
\prref{3.10}(2)). We will see that - given some finiteness condition of the
ring morphism $B\to A$ - this comes down to $A$ being flat as a left $B$-module.\\

Let $V$ be a simple subcomodule of a right $H$-comodule $V$. The sum $M_V$ of all
the subcomodules of $M$ isomorphic to $V$ will be called the {\sl $H$-isotypic
component} of $M$. This sum is a direct sum, and $M_V$ is a semisimple
subcomodule of $M$.\\

We want to describe the $H$-isotypic components of ${}_B\HOM(A,M)$. First we
recall the following Lemma (see \cite[2.14]{19} in the case where $H$ is 
cocommutative). Also recall from \cite[Prop. 2.4.13]{11}
that a simple $H$-subcomodule of an $H$-comodule is finite dimensional.

\begin{lemma}\lelabel{3.16} 
Let $k$ be algebraically closed. Take
$N\in \Mm^H$ and $V\in \Mm^H$
a simple $H$-subcomodule. Then 
$$\Hom^H(V , N)\otimes V\cong N_V$$
as $H$-comodules, and 
$$N=\bigoplus \{N_V~|~V\subset N~{\rm is~a~simple~
subcomodule}\}.$$
\end{lemma}

\begin{remark}\relabel{3.16bis}
If $H$ is commutative, and $N\in {}_A\Mm^H$, then $\Hom(V,N)\cong
V^*\ot N\in {}_A\Mm^H$, so $\Hom^H(V,N)\in {}_B\Mm^H$, and
$\Hom^H(V , N)\otimes V\cong N_V$ is an isomorphism in ${}_B\Mm^H$.
\end{remark}

\begin{lemma}\lelabel{3.17} 
Let $k$ be algebraically closed. Take $M\in
{}_{A}\Mm$ and $V\in \Mm^H$ simple. Then
$${}_B\HOM(A ,M)_V\cong {}_B\HOM(A_{V^*} , M)$$
in $\Mm^H$.
\end{lemma}

\begin{proof} 
Consider the canonical isomorphisms 
\begin{eqnarray*}
&&\hspace*{-2cm}
\Hom^H(V , {}_B\HOM(A, M))
\cong
{}_B\Hom^H(V \otimes A , M)\\
&\cong & {}_B\Hom^H(V \otimes A_{V^*} , M) \cong \Hom^H(V , {}_B\HOM(A_{V^*} , M)).
\end{eqnarray*}
The first and third isomorphism follow from \leref{3.15}(2); the second
follows from the fact that $M$
is a trivial $H$-comodule and from the definition of $A_{V^*}$,
namely, if $W$ is
another simple
$H$-comodule, then $\Hom^H(V\otimes W^* , k)=\Hom^H(V, W^*)=0$ if
$W^* \not= V$. Now it follows from \leref{3.16} that ${}_B\HOM(A , M)_V\cong
{}_B\HOM(A_{V^*} , M)$ as $H$-comodules. If $W$ is another 
simple $H$-comodule, not isomorphic
to $V$, then $\Hom^H(V \otimes A_{V^*} , M)=0$, since $\Hom^H(W \otimes V^* ,
k)=0$. So we find that
$${}_B\HOM(A_{V^*} , M)=\bigoplus_W {}_B\HOM(A_{V^*} , M)_W ={}_B\HOM(A_{V^*} , M)_V.$$
\end{proof}

\begin{remark}\relabel{3.17bis}
If $H$ and $A$ are commutative, then the isomorphism of \leref{3.17} is an isomorphism
in ${}_B\Mm^H$.
\end{remark}

\begin{corollary}\colabel{3.18}  
Assume that $k$ is algebraically closed and that $A$ and $H$ are
commutative. If the functor
$(-)^{{\rm co}H}:\ {}_{A}\Mm^H\to {}_B\Mm$ preserves injectives,
then for every injective left $B$-module $I$, and for every simple
$H$-comodule $V$, ${}_B\HOM(A_{V^*} , I)$ is an injective
left $B$-module.
\end{corollary}

\begin{proof} 
Set $W=V^*$ and $E={}_B\HOM(A_{V^*} , I)$. Then by \prref{3.3}(2), $E$ is an
injective object of ${}_{A}\Mm^H$. Then $M=A \otimes V$ is finitely
generated in
${}_{A}\Mm^H$ and $A$-free. By \reref{3.15bis}, $E
\otimes W\cong
{}_A\Hom(M , E)$ in ${}_A\Mm^H$. By \leref{2.2}, ${}_A\HOM(M
, E)={}_A\Hom(M , E)$, so it follows from \prref{2.14}
that ${}_A\Hom^H(-  , E\otimes W)={}_A\Hom^H(-\otimes_AM , E)$. 
Now ${}_A\Hom^H(- , E):\ {}_{A}\Mm^H\to \Mm$ is exact, so
$E\otimes
W$ is an injective object of ${}_{A}\Mm^H$, and it follows from the
hypotheses that
$(E\otimes W)^{{\rm co}H}$ is injective in ${}_{B}\Mm$. By \leref{3.15}
and \reref{3.15bis}, $(E\otimes
W)^{{\rm co}H}\cong \Hom^H(V , E)$  as a left $B$-module. 
Now, $E_V\cong \Hom^H(V , E)\otimes V$ in ${}_B\Mm^H$, so $E_V$
is injective in ${}_{B}\Mm$. By \leref{3.16}, $E_V={}_B\HOM(A ,
I)_V={}_B\HOM(A_{V^*} , I)$.
\end{proof}

Under the assumptions of \coref{3.18}, if 
$A_{V^*}$ is finitely generated as a left $B$-module, then
${}_B\Hom(A_{V^*} , I)={}_B\HOM(A_{V^*} , I)$ is  an injective left
$B$-module. We will apply this result in \thref{3.19}.

\begin{theorem}\thlabel{3.19} 
Let $k$ be an algebraically closed field and take $A$ and $H$ 
commutative. Assume that the functor
$(-)^{{\rm co}H}:\ {}_{A}\Mm^H\to {}_B\Mm$ preserves injectives,
and let $V$ be a simple $H$-comodule. If $A_V$
is finitely generated as a $B$-module, then $A_V$ is  flat as a left $B$-module. If $A_W$
is finitely generated as a $B$-module, for every simple $H$-comodule $W$, then $A$ is $B$-flat.
\end{theorem}

\begin{proof} Let $I\in {}_{B}\Mm$ be injective. Then by \cite[Prop. 6.5.1]{10}, we have the duality
isomorphism
$${}_B\Hom(\Tor_1^B(M , A_V), I)={}_B\Ext^1(M , {}_B\Hom(A_V, I)).$$
\leref{2.2} and \coref{3.18} show that ${}_BHom(\Tor_1^B(M , A_V), I)=0$ for all $M\in
{}_{B}\Mm$, and it follows that $\Tor_1^B(M , A_V)=0$: it suffices to take
$I={}_BE(\Tor_1^B(M , A_V))$.
\end{proof}

\begin{lemma}\lelabel{3.20} 
Let $A$ be flat as a right $B$-module. Then the functor
$(-)^{{\rm co}H}:\ {}_{A}\Mm^H\to {}_B\Mm$ preserves injectives.
\end{lemma}

\begin{proof} 
We have that ${}_A\Hom^H(A\otimes_B(-) , I)={}_B\Hom(- ,
I^{{\rm co}H})$ in ${}_{B}\Mm$, by \leref{3.2}(2). The
functor $A\otimes_B (-):\ {}_{B}\Mm\to {}_{A}\Mm^H$ is exact. If
$I\in {}_{A}\Mm^H$ is injective, then the functor ${}_A\Hom^H(- , I)$
is exact.
\end{proof}

\begin{lemma}\lelabel{3.21} 
Let $M$ be a finitely generated left $B$-module, and $N\in
{}_{A}\Mm^H$. Then for every $i$, ${}_B\Ext^i(M , N)$ is an
$H$-comodule and
$${}_B\Ext^i(M , N)^{{\rm co}H}={}_B\Ext^i(M , N^{{\rm co}H}).$$
\end{lemma}

\begin{proof} 
Let $\{F_i\}$ be a finitely generated free resolution of $M$.
We can regard each $F_i$ and $M$ as
objects of $_{B}\Mm^H$, with trivial $H$-coaction. It follows
from \cite[Lemma 1.1]{9} that each
${}_B\Hom(F_i , N)$ is an $H$-comodule. Therefore, each ${}_B\Ext^i(M , N)$ is an
$H$-comodule.  Applying \cite[Lemma 1.1]{9} again,
we find that ${}_B\Hom(F_i ,  N)^{{\rm co}H}$=${}_B\Hom(F_i , N^{{\rm co}H})$ for all $i$. The last
assertion follows from
the fact that the functor $(-)^{{\rm co}H}$ commutes with homology.
\end{proof}

\begin{proposition}\prlabel{3.22} 
Let $A$ be finitely generated as a left $B$-module, 
$M\in {}_{B}\Mm$, and $N\in {}_{A}\Mm^H$. Then we have a
spectral sequence
$${{}_A\Ext^H}^i(N , {}_B\Ext^j(A , M))~~ \Rightarrow~~ {}_B\Ext^{i+j}(N^{{\rm co}H} ,
M);~~ i, j \geq 0.$$
If $A$ is left noetherian and $N$ is finitely generated as a left
$A$-module, then
$${}_A\Ext^i(N , {}_B\Ext^j(A , M))^{{\rm co}H}~~ \Rightarrow
~~{}_B\Ext^{i+j}(N^{{\rm co}H} , M);~~ i, j \geq 0.$$
\end{proposition}

\begin{proof} Since $A$ is a finitely generated $B$-module, we have
${}_B\Hom(A , M)={}_B\HOM(A , M)$ and the
first assertion follows from \prref{3.3} and the Grothendieck spectral
sequence for composite functors.
The  second assertion then follows from \coref{2.12}.
\end{proof}

\end{document}